\documentclass[12pt,draftcls,onecolumn]{IEEEtran}

\newtheorem{example}{Example}

\usepackage{color}

\usepackage{amsfonts}
\usepackage{graphicx}      
\newtheorem{theorem}{Theorem}
\newtheorem{lemma}{Lemma}{\bf}{}
{\bf}{}
{\bf}{}
\newtheorem{cor}{Corollary}{\bf}{}
{\bf}{}
\newtheorem{definition}{Definition}{\bf}{}
\newtheorem{assumption}{Assumption}{\bf}{}
{\bf}{}

\begin{document}

\title{Quantized Dissensus in  Networks of Agents subject to Death and Duplication\thanks{Research was supported by grant MURST-PRIN 2007ZMZK5T}}

\author{D. Bauso,\thanks{D. Bauso is with DINFO, Universit\`a di Palermo, 90128 Palermo, Italy,   E-Mail: dario.bauso@unipa.it }
L. Giarr\'e, \thanks{Corresponding author L. Giarr\'e is with
DIAS, Universit\`a di Palermo, 90128 Palermo, Italy, E-Mail:
giarre@unipa.it}  and R. Pesenti\thanks{R. Pesenti is with DMA,
Universit\`a ``Ca' Foscari'' Venezia, 30123 Venezia, Italy,   E-Mail:
pesenti@unive.it } }

 \markboth{IEEE Transactions On AUTOMATIC CONTROL, Vol. XX, No. Y, Month, 200X} {IEEE Transactions On AUTOMATIC CONTROL, Vol. XX, No. Y, Month, 2009}

\maketitle

\begin{abstract}
\emph{Dissensus} is a modeling framework for networks of dynamic agents in competition for scarce resources. 
Originally inspired by biological cells behaviors, it fits also marketing, finance and many other application
areas. Competition is often unstable in the sense that strong agents, those having access to large resources, 
gain more and more resources at the expense of weak agents. Thus, strong agents duplicate when 
reaching a critical amount of resources, whereas weak agents die when loosing all their resources. 
To capture all these phenomena we introduce systems with a discrete time gossip and unstable state dynamics 
interrupted by discrete events affecting the network topology. Invariancy of states
and topologies and network connectivity are explored.   

\end{abstract}

\IEEEpeerreviewmaketitle

\section{Introduction}
Recently, a great interest has been devoted to consensus problems
(see, e.g., \cite{olfati}, \cite{beard} and the literature cited
within). Consensus refers to a scenario where 
a number of agents with interdependent dynamics converge to a common value.
When this happens we say that the agents reach an agreement or \emph{consensus}.
If state $x_i$ is the resource owned by agent $i$, consensus means that 
strong agents, those with large resources, transfer part of these resources to
weak agents. Interdependency has a ``local'' flavor, in the sense that 
each agent can exchange resources only with a subset (say it
\emph{neighborhood}) of adjacent agents. 
In its simplest form the consensus problem can be modeled as an autonomous cooperative
linear system with a Laplacian  dynamic matrix.

In this work, we study a discrete event system (\cite{L03},~\cite{JOAO}) 
that behaves in a complementary manner with respect to the above system. 
We model a competitive scenario where strong agents, 
gain more and more resources at the expenses of weak ones. 
This unstable dynamics can be captured by assuming that 
between two consecutive events, the system evolves as a discrete time autonomous linear system as before, but now 
the dynamic matrix is the opposite of a time dependent Laplacian matrix.
Deviation between neighbors' states/resources increases and for this reason the system is named \emph{competitive} system.

With this premise, it makes sense to assume that when an agent looses all its resources, its
state goes to zero, then it ``dies''.  
The dead agent is then removed from the system together with all its connections. 
\emph{Death} introduces a discrete event into the system dynamics.
  
On the contrary, when an agent gains a critical amount of resources,
then it \emph{duplicates}. The agent, say it \emph{parent}, divides in two new
agents, say it \emph{childrens}. Both the childrens' initial states and
neighborhoods are functions of the parent state and neighborhood
respectively. \emph{Duplication} introduces furtherly a second discrete event into the system dynamics.
Such a modeling framework which we formalize in this paper for the first time, goes under the name of \emph{Dissensus}.
The law ruling the discrete time evolution will be called \emph{dissensus protocol}.

Dissensus is motivated by applications where multiple agents compete for scarce resources. 
In marketing, as an example, agents represent firms or companies and resources are market shares or 
number of customers. Competition often results in firms with a large customer base attracting more and more customers at the expenses of the smaller firms~\cite{VDB2001,Hill2006}. 

Death describes a firm abandoning the market for lack of customers. Duplication captures 
the situation where a firm expands and then divides into smaller firms. 

Duplication may be forced by different causes. For example,
the intervention of an anti-trust authority; the presence of physical constraints in large infrastructures such as airports; the departure of apprentices/young partners from the mother firm to open a new one, in case of service businesses or handicraft industry. 

Other applications are in biological networks, e.g., adjacent cells compete for nutrients \cite{bio,bio1}. The state of an agent describes its size, that in turn is (approximately) proportional to the amount of nutrients the cell subtracts to its adjacent cells. 

Finally, as a theoretical example, let us introduce what we call the \emph{dynamic pie diagram}. In a dynamic pie diagram, the amplitude of each slice evolves over time. The larger slices become ampler and ampler at the expenses of their adjacent slices. Whenever the amplitude of a slice becomes zero, the slice disappears (dies). On the other hand, whenever the amplitude of a slice reaches a given angle the slice is partitioned (duplicates). 

Resources flows are quantized and involve at each time a single pair of agents selected according to a gossip 
rule. Dissensus in continuous time and continuous/unquantized flows have been introduced by the same 
author in~\cite{BGP09}. 

Quantization in consensus problems can be found in~\cite{KBS07,asu}. 
Gossip algorithms in consensus problems have been introduced in~\cite{boyd,fagnani2}. 

The paper is organized as follows. In Section~\ref{sec:CP}, we present the dissensus protocol. In Section~\ref{sec:BasicProperties}, we discuss invariancy of the states and network connectivity.
In Section~\ref{sec:division}, we comment invariancy of topologies under specific death and duplication rules. 
Some final comments are drawn in section~\ref{sec:conclusions}.

\section{The Dissensus Protocol}\label{sec:CP}

The system evolves according to a quantized discrete time dynamics interrupted at \emph{critical times} by some discrete events which we will call \emph{critical events}.
Between each two consecutive critical times~$t_k$ and~$t_{k+1}$, 
with~$t_0 = 0$, the system is made of a set of agents 
$$\label{gamma}\Gamma(t_k) = \{1, \ldots, n(t_k)\},$$ 
whose cardinality~$n(t_k)$  is a function of the time instants~$t_k$.
The agents~$i \in \Gamma(t_k)$ are characterized by states~$x_i$ that assume only discrete values in~$\mathbb{N}$. 
The same agents are organized in a single component \emph{connection network}~$G(t_k)= (\Gamma(t_k),E(t_k))$ whose edgeset~$E(t_k)$ includes all the non oriented pairs~$(i,j) \in \Gamma(t_k) \times \Gamma(t_k)$ of agents 
that bilaterally exchange information about their respective states. 
At the occurrence of critical event both the number of agents and the topology of the connection network are modified.

We describe how the system state $x=\{x_i; i \in \Gamma(t_k)\}$ evolves in the following subsection, whereas
we describe how the system structure is modified at the critical events in Subsection~\ref{subsec:CritEv}. 

Hereafter, for each~$i \in \Gamma(t_k)$,
we call the set  $N_i(t_k)=\{j:(i,j) \in E(t_k)\}$  \textit{neighborhood} of agent~$i$. In addition, for the easy of notation, we omit the dependence on~$t_k$ when there is no risk of ambiguity. Then, e.g., we write~$n$ instead of~$n(t_k)$.

\subsection{Quantized discrete time dynamics}\label{subsec:QuantDyn}

The evolution of each agent state depends on local information meaning that it depends only on the state of the agent' neighbors.
It is governed by a \emph{Quantized Gossip Algorithm} (QGA) in the spirit of what described \cite{KBS07} for cooperative systems. 

\begin{figure}[h!]
\centering
    \includegraphics[scale=0.4]{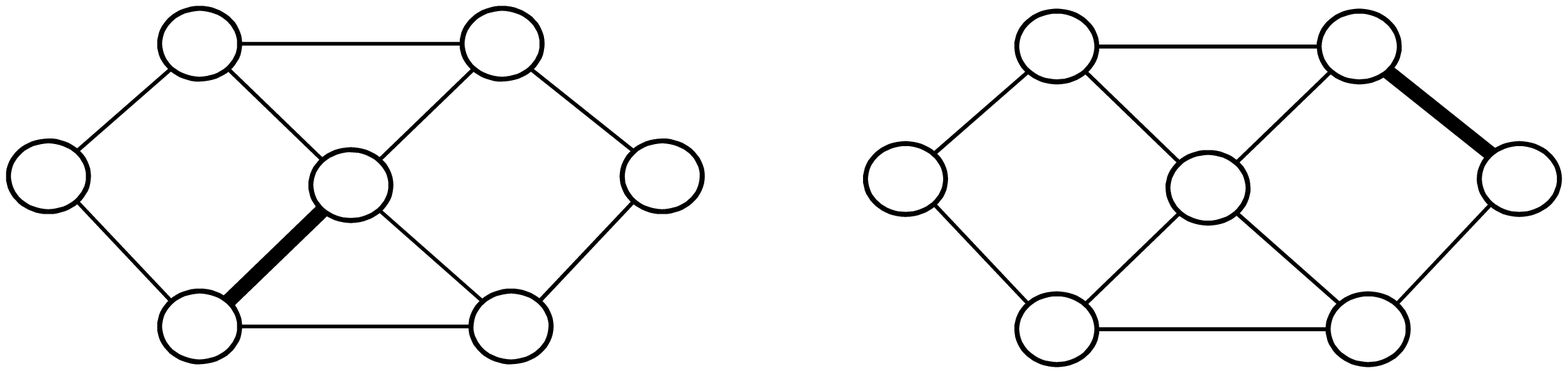}
    \caption{Random selection of links (thick edges) at generic time $r$ (left) and $r+1$ (right).}
   \label{fig:Gossiponetwo}
\end{figure}

At each time~$r$, for any $t_k \leq r < t_{k+1}$,  one edge~$(i,j) \in E$ is
selected, see, e.g., in Fig.~\ref{fig:Gossiponetwo}, and the states~$x_i$ and~$x_j$ of agents~$i$ and~$j$ are updated as follows:

\begin{equation}
x_{i}(r+1) = \left\{ 
\begin{array}{ll}
x_i(r) & \mbox{if $x_i(r) = x_j(r)$} \\
x_{i}(r) + \delta_{ij}(r) & \mbox{if $x_i(r) > x_j(r)$} \\
x_{i}(r) - \delta_{ij}(r) & \mbox{if $x_i(r) < x_j(r)$} 
\end{array} 
\right., 
~ 
x_{j}(r+1) = \left\{ 
\begin{array}{ll}
x_j(r) & \mbox{if $x_i(r) = x_j(r)$} \\
x_{j}(r) - \delta_{ij}(r) & \mbox{if $x_i(r) > x_j(r)$} \\
x_{j}(r) + \delta_{ij}(r) & \mbox{if $x_i(r) < x_j(r)$} 
\end{array} 
\right.,
\label{eq:Goss1}
\end{equation}
for some integer~$\delta_{ij}(r)$ such that $1 \leq \delta_{ij}(r) \leq \min\{x_{i}(r),x_{j}(r)\}$.

Dynamics (\ref{eq:Goss1}) makes the deviation between neighbor agents \emph{diverge} meaning that if, e.g., $x_{i}(r) > x_{j}(r)$, the greater state tends to increase $0 \leq x_{i}(r) < x_{i}(r+1)$, and lower state tends to decrease $x_{j}(r) > x_{j}(r+1) \geq 0$. Also, dynamics (\ref{eq:Goss1}) is \emph{zero-sum}  meaning that a same quantity is gained by $x_{i}$ and lost by $x_{j}$, thus resulting in the local invariance property $x_i(r+1)+x_j(r+1) = x_i(r)+x_j(r)$.
\vspace*{-0.5cm}
Differently from  \cite{KBS07},
 we additionally require that edges are selected in a deterministic way according to the next assumption.
\begin{assumption}\label{as:T}
There exists a big enough integer value~$T$ such that, within each time interval~$T$, either each edge~$(i,j) \in E$ is selected at least once or a critical event occurs. 
\end{assumption}

Had we assumed, in line with \cite{KBS07},  that each edge in~$E$ is selected randomly and with a positive probability, some of the results in this paper would have held in probability.

\subsection{Critical events: death and duplication}\label{subsec:CritEv}

Upper and  lower thresholds are given: the upper threshold is the positive integer value~$B$, the lower threshold is the value zero. A \emph{critical event} occurs at time $t_{k+1}= r$ if an agent~$i$ reaches a threshold at time~$r-1$, i.e, when either $x_{i}(r-1) = 0$ or $x_{i}(r-1) = B$. More formally,
$$t_{k+1}= arg \min\{r > t_k:~\exists i \in \Gamma ~s.t.~x_i(r-1)= 0 ~\vee ~x_i(r-1)= B\}.$$
Hereafter, when it is not stated otherwise, we assume without loss of generality that, within each interval $t_k \leq t < t_{k+1}$, the agents are renumbered so that the agent reaching a threshold value is always agent~$n$. 

On the occurrence of a critical event, the number of agents and the connection network are modified. 
In this context, we denote as~$\Lambda$ (respectively~$\mathcal{E}$) the set of agents (respectively edges) involved in a critical event and still present in the network at the end of the event. 
At~$t_{k+1}$ the system structure modifies as follows. 

\begin{itemize}
	\item If $x_n(t_{k+1 }-1) = 0$, we say that agent~$n$ \emph{dies} and is consequently removed from the system. 
The agents in the neighborhood of~$n$ inherit the connections of~$n$ as depicted in Fig.~\ref{fig:NetworkDeath}. 
\begin{figure}[h!]
\centering
    \includegraphics[scale=0.4]{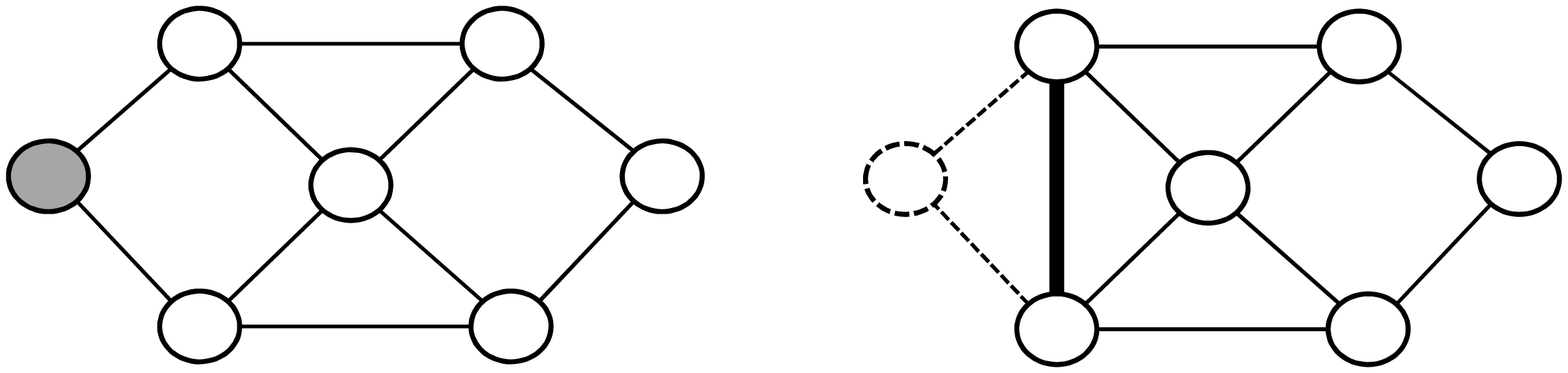}
    \caption{Death of an agent (gray node - left): a new link (thick edge) connects its neighbors (right).}
   \label{fig:NetworkDeath}
\end{figure}
\vspace*{-0.5cm}
In particular, if, for 
each agent $j \in N_n(t_{k})$, we denote with 
$\Lambda_j \subset N_n(t_{k}) \setminus\{j\}$, the set of new neighbors that $j$ inherits from the dead agent $n$, we can write: 
\begin{equation}
\begin{array}{rcl}
		N_j(t_{k+1}) & = & \left\{  \begin{array}{l}
	  (N_j(t_{k})\setminus \{n\}) \cup \Lambda_j  \quad \forall j \in  N_n(t_{k}) \\
	  N_j(t_{k}) \quad otherwise
\end{array} \right., \\
	x_j(t_{k+1}) & = &  x_j(t_{k+1}-1) \quad \forall j \in  \Gamma \setminus \{n\}. 
\end{array}\label{eq:death} 
\end{equation}	

Then, in the occurrence of a death,  the set of involved agents is $\Lambda =  N_n(t_{k})$ and the set of inherited links is ${\mathcal E}=\{(i,j): i \in N_n(t_{k}) \wedge j \in \Lambda_i\}$.

As it will be clearer in Section \ref{sec:BasicProperties}, in order to preserve the connectivity of the network, we require that i) all neighbors of the dead agent are linked to at least one other neighbor, ii) the set of neighbors of the dead agent, together with the inherited links, form a connected subnetwork.
Formally,
    \begin{eqnarray}\label{conn1}
      \mbox{i):  } \bigcup_{j \in N_n(t_{k})} \Lambda_j = \Lambda, \quad \quad \mbox{ii):  }(\Lambda,\mathcal{E}) \mbox{ is connected.}
    \end{eqnarray}
Regarding the connection network~$G(t_{k+1})$, the dead agent~$n$ is removed from the set of agents~$\Gamma(t_{k})$, all links of the dead agent are removed from the set of edges~$E(t_{k})$, finally the new links in~$\mathcal{E}$ are added to~$E(t_{k})$:
    \begin{eqnarray}
    \Gamma(t_{k+1}) = \Gamma(t_{k}) \setminus \{n\}, \quad \quad E(t_{k+1})  &=&  (E(t_{k}) \setminus \{(j,n): j\in N_{n}(t_{k}) \})\cup \mathcal{E}.
    \end{eqnarray}

  \item If $x_n(t_{k+1}-1) = B$, the parent agent~$n$ \emph{\emph{divides}} producing two children agents~$n$ and~$n+1$. The two children agents inherit the parent connections and state (see, e.g., Fig.~\ref{fig:Death1Dupli2}). 
\begin{figure}[h!]
\centering
    \includegraphics[scale=0.4]{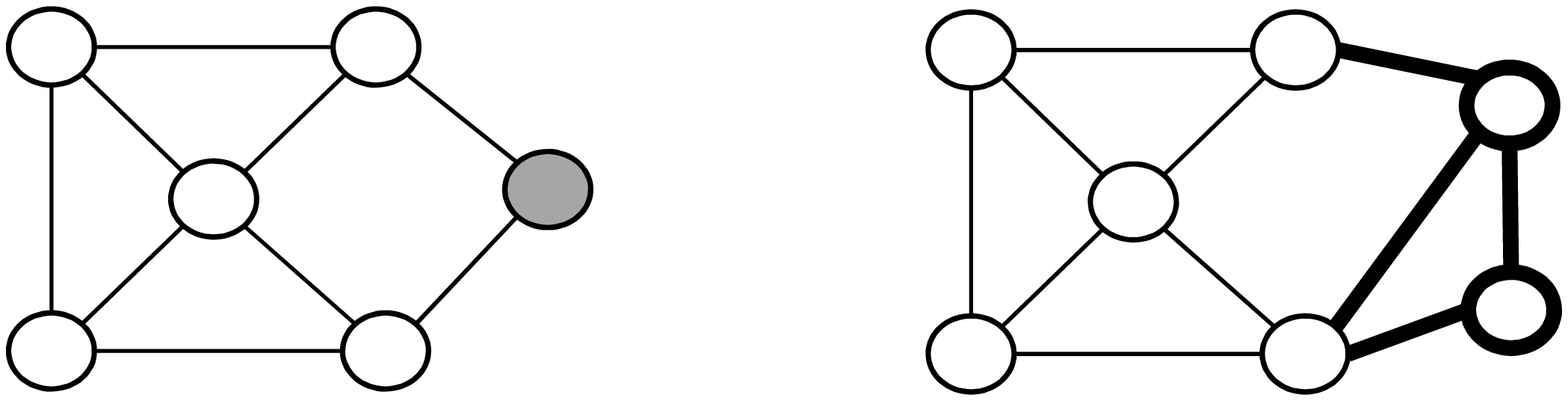}
    \caption{Duplication of one agent (gray node - left): new agents and links are added (thick edges and nodes - right).}
   \label{fig:Death1Dupli2}
\end{figure}
\vspace*{-0.5cm}
This means that, if we denote with 
$\Lambda_{n}$ and~$\Lambda_{n+1}  \subseteq N_n(t_{k}) \cup \{n, n+1\}$ the new neighbors of the  two children agents~$n$ and~$n+1$, we can write:
\begin{equation}
\begin{array}{rcl}
		N_j(t_{k+1})  &= & \left\{  \begin{array}{ll}
	\Lambda_n   & \mbox{for $j =n$}\\
	\Lambda_{n+1} & \mbox{for $j =n+1$} \\
	N_j(t_{k}) \cup \{n\} & \mbox{if $j \in \Lambda_n$, $j \not= n+1$}\\
	(N_j(t_{k})\setminus \{n\}) \cup \{n+1\} & \mbox{if $j \in \Lambda_{n+1}$, $j \not= n$}\\
	N_j(t_{k}) & \mbox{otherwise}
\end{array} \right., \\
	x_j(t_{k+1})  &= &  \left\{  \begin{array}{ll} 
	\alpha & \mbox{for $j =n$}\\
	\beta & \mbox{for $j =n+1$} \\
  x_j(t_{k+1}-1) & \mbox{otherwise}
	\end{array} \right.. 
\end{array}
\label{eq:division}
\end{equation}	
Then, in the occurrence of a duplication, the set of involved agents is $\Lambda =  N_n(t_{k}) \cup \{n, n+1\}$ and the set of the links to the  children agents is $\mathcal{E} = \{(n,i): i \in  \Lambda_{n}\} \cup \{(n+1,i): i \in  \Lambda_{n+1}\}$.

In order to preserve the connectivity of the network, we impose that one of the following two rules hold
\begin{equation}\label{connectivity1} \label{connectivity2}
\begin{array}{rcl}
(\Lambda_{n}\cup \Lambda_{n+1} = \Lambda) \quad \mbox{ or }   \quad
(\Lambda_{n}\cup \Lambda_{n+1} = N_n(t_{k}),~\Lambda_{n}\cap \Lambda_{n+1} \not = \emptyset). 
\end{array}
\end{equation}

Furthermore, in order to maintain the sum of the states invariant, we impose also that the values~$\alpha$ and~$\beta$ are integers such that $\alpha \geq \beta >0$ and  \begin{equation}\label{inv}\alpha+\beta = B = x_n(t_{k}).\end{equation} 

Regarding the connection network~$G(t_{k+1})$, a new agent~$n$, replacing the old one, and an additional agent~$n+1$ are added to the set of agents~$\Gamma(t_{k})$, all the links to the father agent~$n$ are removed from the set of edges~$E(t_{k})$ and the set~$\mathcal{E}$ of new links to the children agents are added to~$E(t_{k})$:
		\begin{eqnarray}\label{dupli111}
		\Gamma(t_{k+1})  =  \Gamma(t_{k}) \cup \{n+1\}, \quad 	E(t_{k+1}) & = & (E(t_{k}) \setminus \{(j,n): j\in N_{n}(t_{k}) \} )\cup \mathcal{E}.
	  \end{eqnarray}

\end{itemize}

\begin{figure}[h!]
\centering
    \includegraphics[scale=0.40]{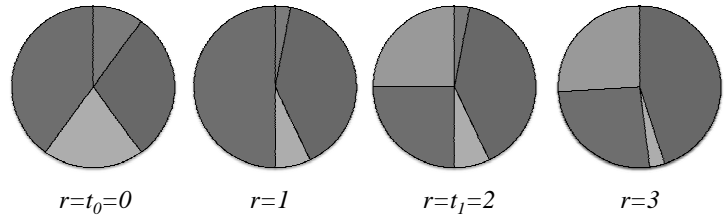}
    \caption{The evolution of a dynamic pie diagram.}
   \label{fig:PieDiagram}
\end{figure}

Figure~\ref{fig:PieDiagram} shows the possible evolution of a dynamic pie diagram. An agent/slice divides at time~$r=2$, having reached an amplitude of $B=180^o$ at the previous time instant. Differently, an agent/slice dies at time $r=3$.

\vspace*{-0.5cm}
In the next session, we discuss invariancy and connectivity associated to death  and duplication rules.

\section{Connectivity and Invariancy}\label{sec:BasicProperties}

In this section, we introduce some bounds on the number~$n(t_k)$ of agents. 
To this end, we initially discuss the  connectivity and invariancy properties that are exploited in the bounds determination.

Preserving the agents' connectivity motivates the choice for death~(\ref{conn1}) and duplication rules~(\ref{connectivity1}). 
Before discussing them in details, we need the following definition.

\begin{definition}(\emph{Local connectivity})
$G(t_{k+1})$ is \emph{locally connected} if the subnetwork $(\Lambda,\mathcal{E})$ 
is connected. 
\end{definition}
Note that we can impose the local connectivity of~$G(t_{k+1})$ exploiting the knowledge of~$N_n(t_{k+1})$, i.e., the information available to the agent~$n$ just before its death or duplication. 
Conditions for the local connectivity to imply the connectivity of~$G(t_{k+1})$ are established in the next theorem.

\begin{theorem} 
If~$G(t_k)$ is connected, then both the
death~(\ref{conn1}) and duplication rules~(\ref{connectivity1})
are sufficient conditions for the connectivity of~$G(t_{k+1})$. 
These rules become even necessary to guarantee the connectivity of~$G(t_{k+1})$ 
on the basis of the only knowledge of neighborhood~$N_n(t_k)$ available to agent $n$ and disregarding any information
on other  agents' neighborhoods~$N_i(t_{k})$ for~$i \in \Gamma(t_k)\setminus \{n\}$.
\end{theorem}

\begin{proof}
The proof is based on the following two facts:
\begin{itemize}
	\item $G(t_{k+1})$ differs from~$G(t_{k})$ only in the subnetwork induced by the agents in~$\Lambda$,  
  \item both the death~(\ref{conn1}) and duplication rules~(\ref{connectivity1})
guarantee the local connectivity of~$G(t_{k+1})$.   
\end{itemize}
This latter fact hold true by definition of rule~(\ref{connectivity2}) in case of a death. 
In the case of a duplication, we note that~$\Lambda = N_{n}(t_k) \cup \{n, n+1\}$ and that
both rules~(\ref{connectivity1}) imply that each~$i \in N_{n}(t_k)$ is connected with either~$n$ or~$n+1$. 
In addition, if the first condition of~(\ref{connectivity1}) holds, then $n+1\in \Lambda_n$ and $n\in \Lambda_{n+1}$, hence, $\mathcal{E}$ includes~$(n,n+1)$.
If the second condition of~(\ref{connectivity2}) holds, then there exists~$i \in N_{n}(t_k)$ such that ~$\mathcal{E}$ includes~$(n,i)$ and~$(n+1,i)$. In both cases, we can conclude that each~$i \in \Lambda$ is connected to~$n$ through the edges in~$\mathcal{E}$. By transitivity $(\Lambda, \mathcal{E})$ is connected.

(Sufficiency) As~$G(t_k)$ has a single component, we know that, in~$G(t_k)$,  each~$i \in \Gamma(t_k) \setminus (N_n(t_k) \cup \{n\})$ is connected to at least an agent $j_i \in N_n(t_k)$ through a path that does not include agent~$n$ and any other agent in~$N_n(t_k)$ different from~$j_i$. These paths connecting~$i$ with~$j_i$ remain identical in~$G(t_{k+1})$. Then, as~$N_n(t_k) \subseteq \Lambda $, all the agents $i \in \Gamma(t_k) \setminus (N_n(t_k) \cup \{n\})$ are connected  in~$G(t_{k+1})$ with at least an agent in~$\Lambda$. As~$(\Lambda,\mathcal{E})$ is connected, all the agents in~$\Lambda$, and hence the agents in $i \in \Gamma(t_k) \setminus (N_n(t_k) \cup \{n\})$, are connected to each others in~$G(t_{k+1})$. Finally, as the set of the agents in~$G(t_{k+1})$ is $\Gamma(t_{k+1}) = (\Gamma(t_k) \setminus (N_n(t_k) \cup \{n\})) \cup \Lambda$, we can conclude that~$G(t_{k+1})$ is connected.

(Necessity) If we have no knowledge of neighborhoods~$N_i(t_{k})$ for $i \in \Gamma({t_k})\setminus \{n\}$, we must conservatively assume that all the links between~$n$ and the agents~$N_n(t_{k})$ are necessary for the connectivity of~$G(t_{k})$, since~$G(t_{k})$ may be, e.g., a tree. Consequently, to guarantee the connectivity of~$G(t_{k+1})$, we must guarantee that the edges in $\mathcal{E}$ keep the agent in~$N_n(t_{k})$ connected to each others and to the new agents eventually added in $t_{k+1}$. In other words, we require  that~$(\Lambda,\mathcal{E})$ is connected.

\end{proof}

\begin{figure}[h!]
\centering
\begin{minipage}[c]{.40\textwidth}
\centering
\includegraphics[scale=0.4]{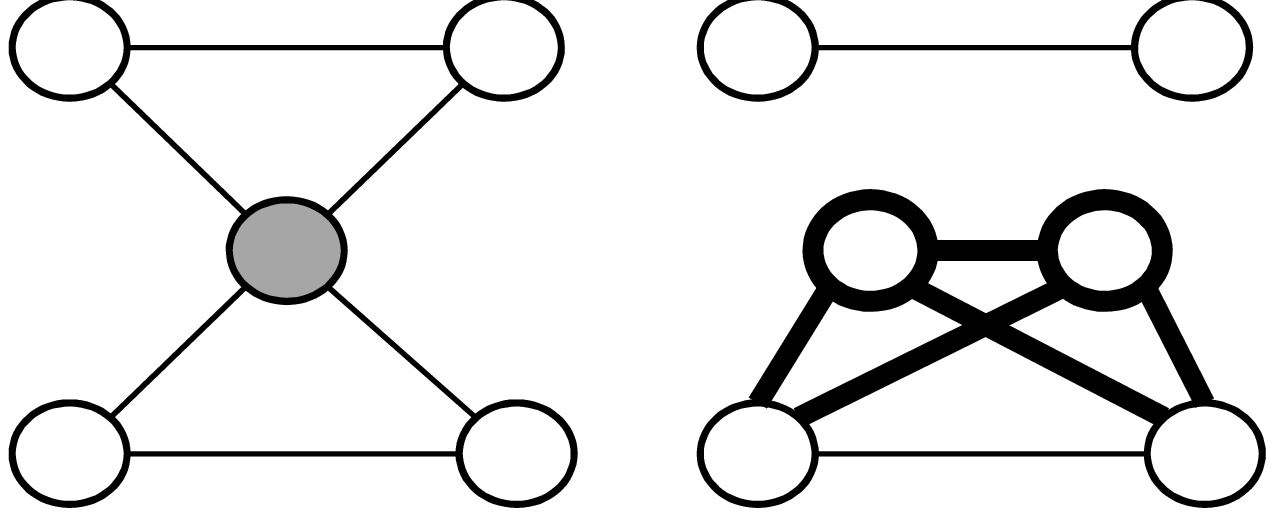}
\end{minipage}%
\hspace{10mm}%
\begin{minipage}[c]{.40\textwidth}
\centering
\includegraphics[scale=0.4]{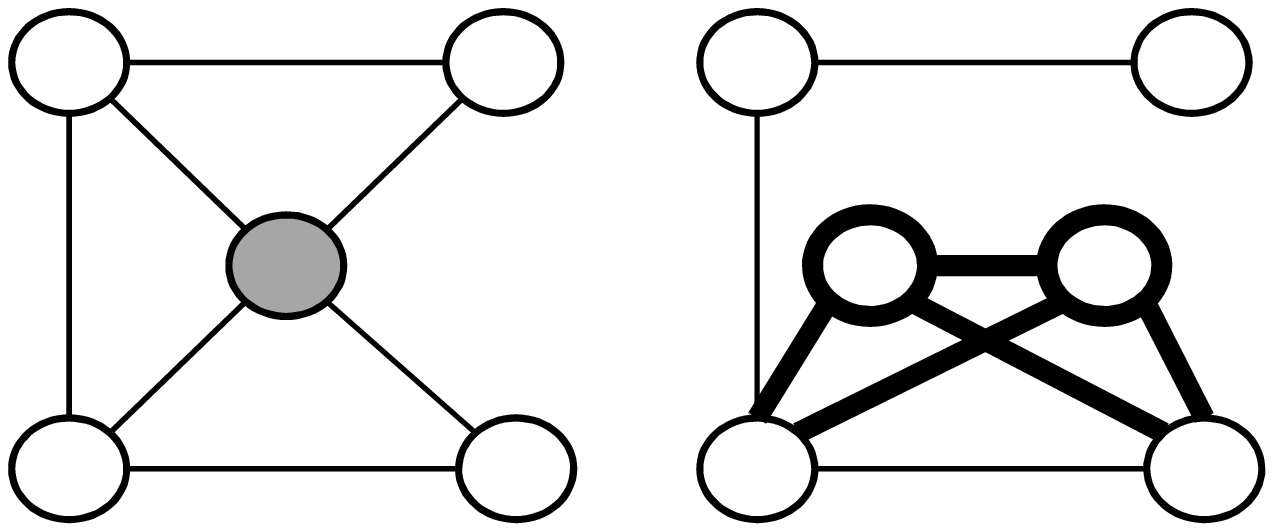}
\end{minipage}
\caption{Different duplication events. Duplicating agent in gray, children agents in thick line. Cases in which rules~(\ref{connectivity1}) are  necessary (left networks) and not necessary to the connectivity of~$G(t_{k+1})$ (right networks).}
   \label{fig:examples}
\end{figure}

Figure~\ref{fig:examples} left shows that if~(\ref{connectivity1}) does not hold~$G(t_{k+1})$ may turn to be disconnected. 
However, in general, rules~(\ref{conn1}) and~(\ref{connectivity1}) are not necessary to the connectivity of~$G(t_{k+1})$.
As an example, in case of a duplication, the connectivity of~$G(t_{k+1})$ holds even if $\Lambda_j = \{j\}$ for all $j \in  N_n(t_{k})$, if all $j \in N_n(t_{k})$ are adjacent to a common agent $v \notin N_n(t_{k})$ (see, e.g., Fig. \ref{fig:examples} right).


As regards invariancy, we observe  that dynamics described in (\ref{eq:Goss1}), (\ref{eq:death}) and (\ref{eq:division})  imply that the sum of the agents' states is a constant value over time, i.e.,  
\begin{equation}
\sum_{i \in \Gamma(r)} x_i(r) = \chi.
\label{eq:StateInvariance}
\end{equation}

Invariancy is a necessary condition to avoid that the number of agents neither diverges to infinity nor converges 
to one. 
Assume that in consequence of a duplication at time $t_k$ the sum of the states decreases over time,
i.e, $\alpha + \beta < x_n(t_k-1)$ or that an agent may die having a state  $x_n(t_k-1) >0$. 
It is immediate to see, that under these hypotheses either the system reaches consensus or $\lim_{k \rightarrow \infty} n(t_k) = 1$. 
Differently, if the sum of the states increases over time, which corresponds to saying $\alpha + \beta > x_n(t_k)$, then either the system reaches consensus or $\lim_{k \rightarrow \infty} n(t_k) = \infty$.

Invariancy is also used in the next result to establish that the number of agents never goes above neither below
certain bounds.

\begin{theorem}\label{th:Bounds}
Let $t_h$ be the critical time, possibly infinite, at which the first duplication occurs. Then the following bounds holds on the number of agents: 
\begin{equation} \left\lceil \frac{\chi}{B} \right\rceil \leq n(t_k) \leq 
\left\{ \begin{array}{ll}
\chi - B +2, &  \forall t_k \geq t_h  \\
\chi, &  \forall t_k < t_h  
\end{array}\right..		
\label{eq:nLowerBound}
\end{equation}
\end{theorem}

\begin{proof}
The lower bound is a straightforward consequence of the invariancy (\ref{eq:StateInvariance}).
As regards to the upper bounds, we can have $n(t_0) = \chi$ in the following trivial situation.
At the initial time $t_0 = 0$ the agents' states have value $x_i(t_0) = 1$ for all $i \in \Gamma(t_0)$.
Differently,  from $t_h$ on, the maximum number of agents  can be reached only when a duplication occurs. In this case, the agents' state invariancy imposes 

\begin{equation}
\chi = \sum_{i \in \Gamma(t_{k})} x_i(t_k) \geq n(t_{k}) - 2  + \alpha + \beta = n(t_{k}) - 2  + B.
\label{eq:LB1Q}
\end{equation}
In (\ref{eq:LB1Q}), the inequality  holds as, after a duplication event, at least one of the system agents has state equal to $\alpha$, a second agent has state equal to $\beta$, and the remaining agents have state $x_i(t_k) \geq 1$. Finally, the bound in (\ref{eq:nLowerBound}) trivially follows from (\ref{eq:LB1Q}).
\end{proof}

Note that the fact that $n(t_k)$, and hence the cardinality of  $E(t_k)$, is  bounded from the above for all $t_k \geq 0$ guarantees the possibility that the value~$T$ required by Assumption~\ref{as:T} exists.

Invariancy and connectivity play a critical role when proving that under a QGA the intervals between consecutive critical events is finite. Indeed, this fact holds true unless the agents' reach \emph{consensus}.
We say that the agents reach consensus at time $r$ if all states are equal, i.e., $x_i(r) = x_j(r)$ for all $i,j \in \Gamma$. 

The next lemma is a prtelude to show that the interval between consecutive critical events is finite.

\begin{lemma}\label{lem:cv}
Consensus can be reached only at times $t_k$, for $k = 0, 1, \ldots$.
\end{lemma}
 
\begin{proof}
For each time $r = t_k, t_k +1, \ldots, t_{k+1}-1$, let $imin(r)=arg \min_{i \in \Gamma}\{x_i(r)\}$ be the agent with minimum state,  $imax(r)=arg \max_{i \in \Gamma}\{x_i(r)\}$ be the agent with maximum state.   By definition of QGA we have that the values of $x_{imin(r)}(r)$, respectively  $x_{imax(r)}(r)$,  may vary but cannot increase, respectively decrease, for  $t_k \leq r < t_{k+1}$. As a consequence, the difference between $x_{imin(r)}(r)$ and  $x_{imax(r)}(r)$  cannot reduce to zero for  $t_k < r < t_{k+1}$.
\end{proof}

\begin{theorem}\label{th:CritEv}
Assume that the agents do not reach consensus in~$t_k$, then a critical event will occur in a finite time~$t_{k+1}$.
\end{theorem}
\begin{proof}
We prove the statement by contradiction. In particular, we show that if no critical event  occurs in a finite time the system state diverges. To this end, we assume that the agents have not reached consensus at time~$t_k$ and that no critical event occurs from~$t_k$ on.
Let $i$, $j$ be the agents selected by the QGA at the generic time $r > t_k$ and consider difference 
$$
||x(r+1)||^2 - ||x(r)||^2 = \left\{ \begin{array}{ll} 2\delta_{ij}(r)^2 + 2\delta_{ij}(r)(x_{i}(r)-x_{j}(r)) \geq 4 & \mbox{if $|x_i-x_j| \not = 0$}\\
0 & \mbox{otherwise}
\end{array}  \right. .
$$
In view of Assumption~\ref{as:T}, we have that $||x(s)||^2 - ||x(t_k)||^2 \geq 4\left\lfloor \frac{s-t_k }{T}\right\rfloor$   for all~$s\geq t_k$,
as, if i)~no critical events occurs after~$t_k$, ii)~$G(t_k)$ is connected and iii)~the agents have not reached consensus in~$t_k$, within every~$T$ instants at least a couple of agents with different states is selected.
Then, as $\lim_{s \rightarrow \infty}4\left\lfloor \frac{s-t_k }{T}\right\rfloor = \infty$, we have that $\lim_{s \rightarrow \infty} ||x(s)||^2 = \infty$. The value of the latter limit is in contradiction with the fact that $0 \leq x_i(s) \leq B$ must hold for any $s \geq 0$. Hence, either  the agents have reached consensus in~$t_k$ or a critical event occurs after~$t_k$.
\end{proof}

An immediate consequence of Theorems~\ref{th:Bounds} and~\ref{th:CritEv} is that a first duplication event occurs at time $t_h$ after at maximum $\chi-2$ death events. In addition, the $t_h$ is finite unless the agents reach consensus at time $t_0$ or at one the above death events.

In the rest of this section,  we analyze possible scenarios where the agents reach consensus.
\begin{example}
Consider systems evolving according a protocol in which the duplication rule~(\ref{eq:division}) requires that $\alpha = \beta = B/2$ and that both the two children inherit all their parent connections and connect to each other (see, e.g., rule~(\ref{eq:likeFather}) discussed in the next section).  Such systems may reach consensus. Note that in this case each duplication generates twin agents.
Consider, as an example, the evolution of a system where $n(0)=2$ and $x_1(0) = B/2$ and $x_2(0) = 3B/4$. By Theorem~\ref{th:Bounds}, we have $n(t_k) \geq 2$ for all $t_k$. We can easily define a QGA such that at time $t_1$, we have $x_1(t_1) = B/4$ and $x_2(t_1) = x_3(t_1) = B/2$. At time $t_2$, the agent~$1$ dies (we have not renumbered the agents for the easy of exposition) and the remaining two twin agents reach an equilibrium corresponding to a state value equal to $5B/8$.
\end{example}

The next theorem establishes some necessary conditions for the agents to reach consensus.
\begin{theorem}\label{th:cv}
Consensus can be reached provided that  $x_i(t_k) = x_j(t_k) = \frac{\chi}{n(t_k)} < B $ for all $i,j \in \Gamma$.
Furthermore, 
$$\frac{\chi}{n(t_k)} \geq \left\{\begin{array}{ll}
\frac{\chi}{n(t_0)} & \mbox{for $t_0 \leq t_k < t_h$} \\
\alpha	& \mbox{for $t_k \geq t_h$}
\end{array}\right.,$$ 
where $t_h$ is the critical time of the first duplication event.
\end{theorem}
 
\begin{proof}
The first part of the theorem is a trivial consequence of the definition of consensus and of the invariancy property.

The first lower bound is true as, for $t_0 \leq t_k < t_h$,  we have $n(t_k+1) = n(t_k) -1$ and, obviously, the invariancy property holds.
In proving the second lower bound,  we observe that, for $r \geq t_h$, we have $x_{imax(r)}(r)= \max_{i \in \Gamma(t_k)}\{x_i(r)\} \geq \alpha$. Actually, any QGA prevents $x_{imax(r)}(r)$ from decreasing between two consecutive critical duplication events.
\end{proof}

A direct consequence of Theorem~\ref{th:cv} is that there exist initial states~$x(0)$ such that the agents will never reach a consensus as stated in the following corollary.
\begin{cor}
Let $t_h$ be the critical time of the first duplication event.
No consensus can be attained for $t_0 \leq t_k < t_h$, if $\chi$ is not divisible for some integer less than or equal to $n(t_0)$. 
No consensus can be attained for $t_k \geq t_h$, if $\chi$ is not divisible for some integer greater than or equal to $\alpha$. 
\end{cor}
\begin{proof}
The  corollary thesis is an immediate consequence of of Theorem~\ref{th:cv} and the fact that both the number of agents that reach consensus and their states must assume integer values.
\end{proof}

%
We conclude this section by observing that trajectory $x(t)$ may become periodic in the long run when the agent do not reach a consensus. As an example, this situation occurs if the system evolution evolves according to death/duplication rules that depend on the agents' states and on the connection network topology and no stochasticity or time dependency is allowed.Indeed, we can observe that  our systems evolve in a space defined by the number of agents, their states and the topologies of their possible connection networks.As both the number of agents and their states may assume only integer values bounded from above, then 
the number of topologies of the possible connection networks is also finite, hence the thesis follows.

\section{Topology invariancy} \label{sec:division}

Let us now investigate conditions on death and duplication rules that result in the invariancy  of specific network topologies like complete, hole or chain networks. 
As we study the asymptotic evolution of a system, let us initially observe that Theorem~\ref{th:Bounds} implies that there always exists a critical time~$t_q$ an two integer values $\underline{n}$ and $\bar n$ such that, for each $\tilde n$ satisfying  
\begin{equation}
\left\lceil\frac{\chi}{B} \right\rceil\leq \underline{n} \leq \tilde n \leq \bar n \leq \chi,
\label{eq:nbar}
\end{equation}
we have an infinite sequence of critical time instants $\sigma(\tilde n) = \{t^{\tilde n}_k:  t^{\tilde n}_k \geq t_q\}$ such that $n(t^{\tilde n}_k) = \tilde n$.
In this context, we say that the \emph{density of the connection network} of a system does not increase asymptotically if there exists a value~$\hat t \geq 0$ such that for each~$\tilde n$ satisfying condition~(\ref{eq:nbar}), the value $|E(t^{\tilde n}_k)|$ does not increase over the critical times in~$\sigma(\tilde n)$ and greater than or equal to~$\hat t$.

We start by considering the
 duplication rule which equally divides the  parent's connections between the two children. Formally,
\begin{equation}
	\Lambda_n = pick(N_n(t_k)) \cup \{n+1\} \quad \mbox{ and } \quad 
	\Lambda_{n+1}= N_n(t_k) \setminus \Lambda_n \cup \{n\},
\label{eq:Part}
\end{equation}
where function $pick(N_n(t_k))$ returns a random subset of  $\left\lfloor |N_n(t_k)|/2 \right\rfloor$ elements of $N_n(t_k)$.

\begin{figure}[h!]
\centering
    \includegraphics[scale=0.4]{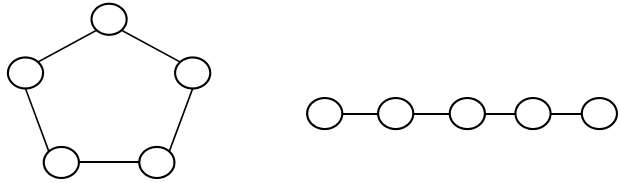}
    \caption{A hole network and a chain network.}
   \label{fig:HoleChain}
\end{figure}

With such a rule, if $G(t_0)$ is a hole network, respectively a chain network (see Fig.~\ref{fig:HoleChain}), then $G(t_k)$ are hole network, respectively chain network, for every $t_k$, whatever death rule is implemented. Again, this corresponds to saying that hole and chain networks are invariant as established in the next theorem. Before giving the theorem, we recall that a network $G(t_k)$ is a hole, respectively a chain network, if it is connected and all the agents have degree two, respectively all the agents have degree two a part from two agents at the extreme of the chain whose degree is one (see, e.g., Fig. \ref{fig:HoleChain}). 

\begin{theorem}
Hole and chain networks are invariant under duplication rule~(\ref{eq:Part}).
\end{theorem}
\begin{proof}
The idea of the proof is that critical events do not change, in general, the degree of the agents. 
This is always true unless the dying or duplicating agent $n$ is an extreme of the chain. Actually,
when an extreme agent dies, its neighbor, whose degree changes from two to one, becomes the new extreme.   
Similarly, when an extreme agent duplicates, one children is the new extreme with degree one, while the 
other children has degree two as all other agents.   
\end{proof}

\vspace*{-0.5cm}
Consider now the death rule~(\ref{eq:death}) that assign all the connections of a 
dead agent to just one of its neighbors. 
This can be formally described by a death rule that, for each $j \in N_n(t_k)$, satisfies: 
\begin{equation}
	\Lambda_{j} =  \left\{ \begin{array}{ll}
	 N_n(t_k) \setminus \{j^*\} & if~j = j^* \\
   \{j^*\} & if~j \not = j^* \\
   \end{array} \right.
	\label{eq:death1}
\end{equation}
where $j^* \in N_n(t_k)$ is arbitrarily  picked.

From a practical point of view, these choices for $\Lambda_{j}$ are the simplest ones to implement that guarantee the connectivity of $G$ after the removal of the dying agent.

\begin{theorem}\label{th:lowden}
Under duplication rule~(\ref{eq:Part}) and death rule~(\ref{eq:death1}), both the density and the number cycles of a connection network cannot increase asymptotically. 
\end{theorem}

\begin{proof}the thesis is trivially true if the agents reach consensus. Otherwise, it is immediately to see that, if at~$t_k$ a  duplication occurs, $|E(t_k)| = |E(t_{k-1})| + 1$; if at~$t_k$ a death occurs, 
\begin{equation}
|E(t_k)| \leq |E(t_{k-1})| - 1
\label{eq:ineq}
\end{equation} as at least the connection~$(n,j^*)$ is not substituted by a new connection. 
Now, observe that, for each~$\tilde n$ satisfying condition~(\ref{eq:nbar}), an equal number of death and duplication events must occur between each two critical times $t^{\tilde n}_k$ and $t^{\tilde n}_s \in \sigma(\tilde n)$, with $t^{\tilde n}_s >t^{\tilde n}_k$, hence $E(t^{\tilde n}_s) \leq  E(t^{\tilde n}_k)$. Then, we can conclude that the density of the connection network cannot increase over time. 
Finally, the observation that the considered  duplication and death rules forbid the creation of new cycles in the connection network concludes the proof.
\end{proof}

Authors' computational experiments  show that the density of the connection network usually indeed decreases until $G$ presents a single or no cycle at all if the assumptions of Theorem~\ref{th:lowden} hold and the agents do not reach consensus. Although no formal proof confirms such observations, it appears reasonable to expect that the the connection network becomes sparser and sparser. Indeed, inequality~(\ref{eq:ineq}) in the proof of Theorem~\ref{th:lowden} holds strictly in the occurrence of a death event whenever $N_n(t_k) \cap N_{j^*}(t_k) \not = \emptyset$, situation quite common if the network is not sparse.

Let us now consider the  duplication rule which makes both children inherit all the parent connections and connect to each other: 
\begin{equation}
\begin{array}{rcl}
		\Lambda_n  =  N_n(t_k) \cup \{n+1\} \quad \quad \mbox{and} \quad \quad
	\Lambda_{n+1}  =  N_n(t_k) \cup \{n\}. 
\end{array}
\label{eq:likeFather}
\end{equation}

Observe that the above rule preserves complete connectivity of the subnetwork induced by the critical set
after a duplication event. This also induces the fact that if $G(t_0)$ is a complete network then $G(t_k)$ are complete networks for every $t_k$, whatever death rule is implemented. This corresponds to saying that the complete network is \emph{invariant} under the duplication rule (\ref{eq:likeFather}) as established in the next theorem.
\begin{theorem}
The complete network is \emph{invariant} under the duplication rule (\ref{eq:likeFather}).
\end{theorem}

\begin{proof}
Assume $G(t_k)$ is a complete network. If a death occurs at time $t_{k+1}$, each agent is already adjacent to all the other agents. If a duplication occurs at time $t_{k+1}$, because of rule~(\ref{eq:likeFather}), the two new agents are adjacent to each other and to all the other agents.
\end{proof}

The authors' computational experiments  show that, for a generic systems with $G(t_0)$ not complete, the density of the connection network usually increases until $G$ becomes a complete network if the rule~(\ref{eq:likeFather}) is applied and the agents do not reach consensus, see \cite{BGP-ve09}.

\section{Final Remarks}\label{sec:conclusions}
In this paper we have introduced systems of competitive agents exchanging quantized flows according to a dissensus protocol. We have discussed some invariancy properties of the systems and observed that the number of agents is bounded from below and above over time. 

As a future direction of research we are oriented toward the definition of more general rules for the inheritance of the connections. Indeed, the main limit of the current version of the dissensus protocol seems the fact that the inheritance rule of the connections applies only when an agent dies. In many real competitive systems the larger agents may acquire the connections of the smaller agents when the latter ones are still alive.  Consider the following two examples. The number of the cells adjacent to (that is, in the neighborhood of) a cell~$i$ of a tissue is a function of the size of the cell~$i$ and not only of the possible occurrence of a death of an adjacent cell.   Similarly, the rules used for the dynamic pie diagram in Fig.~\ref{fig:PieDiagram} cannot trivially be extended to a general dynamic tessellation.



\end{document}